\documentclass[12pt,a4paper,russian,notitlepage]{article}
\usepackage[russian]{babel}
\usepackage{amsmath,amssymb,amsthm}

\newtheorem*{theorem*}{Теорема}
\newtheorem{theorem}{Теорема}
\newtheorem{lemma}{Лемма}
\newtheorem*{lemma*}{Лемма}

\newtheorem*{remark*}{Замечание}

\newcommand{\eqdef}{\mathrel{\mathop=}:}

\begin{document}
\begin{center}{\Large Алгебраическое доказательство линейных теорем Конвея--Гордона--Закса и ван Кампена--Флореса}\footnote{Доказательство возникло в результате обсуждений на семинарах по дискретному анализу на факультете инноваций и высоких технологий Московского физико-технического института.}
\end{center}

\vspace{0.1cm}

\begin{center}{\large  И. И. Богданов, А. Д. Матушкин}
\end{center}

\vspace{0.5cm}

{\bf Abstract.}
In this paper we present short algebraic proofs of the Linear Conway--Gordon--Sachs and the Linear van Kampen--Flores theorems in the spirit of the Radon theorem on convex hulls.

{\bf Theorem.} {\it Take any $n+3$ general position points in $\mathbb{R}^n$. If $n$ is odd, then there are two linked $(n+1)/2$-simplices with the vertices at these points. If $n$ is even, then one can choose two disjoint $(n+2)/2$-tuples such that the interiors $(n/2)$-simplices with the vertices at these $(n+2)/2$-tuples intersect each other.}

This theorem is interesting even in case of small dimensions.

\section{Введение}
\label{intro}
В данной работе предлагаются короткие алгебраические доказательства классических теорем о пересечении и зацепленности.
Полные версии теорем Конвея--Гордона--Закса и ван Кампена--Флореса были сформулированы соответственно для непрерывных вложений полного графа на $n+3$ вершинах в $\mathbb{R}^n$ (точнее, авторы доказали эту теорему для шести точек в $\mathbb{R}^3$) и непрерывных вложений $n$-мерного остова $(2n+2)$-мерного симплекса в $\mathbb{R}^{2n}$ (см.~\cite{Prasolov}). Мы обсуждаем линейные (то есть более частные) аналоги этих теорем. Отметим, что известны другие простые доказательства рассматриваемых теорем (см., например,~\cite{S},~\cite{SS},~\cite{Zimin}), основанные на идее понижения размерности. Теоремы Конвея--Гордона--Закса и ван Кампена--Флореса интересны уже в случае малых размерностей. Мы сначала сформулируем их для трехмерного и четырехмерного пространств соответственно, а затем обобщим формулировки и представим доказательство обобщенных теорем.

\bigskip
Пусть $\Delta$ и $\Delta'$ --- два треугольника в пространстве, среди шести вершин которых никакие 4 не лежат в одной плоскости. Будем говорить, что эти треугольники {\it зацеплены}, если контур треугольника $\Delta$ пересекает внутренность треугольника $\Delta'$ в единственной точке.

\begin{theorem}[Линейная версия теоремы Конвея--Гордона--Закса, частный случай]
\label{CGS}
Для любых 6 точек в трехмерном пространстве, никакие 4 из которых не лежат в одной плоскости, найдутся два зацепленных треугольника с вершинами в этих точках.
\end{theorem}

\begin{theorem}[Линейная версия теоремы ван Кампена--Флореса, частный случай]
\label{vKF}
Среди любых 7 точек в четырехмерном пространстве, никакие 5 из которых не лежат в одном трехмерном аффинном подпространстве, можно выбрать две непересекающиеся тройки точек такие, что треугольники с вершинами в них пересекаются.
\end{theorem}

\begin{remark*}
В оригинальных формулировках теорем ~\ref{CGS} и~\ref{vKF} вместо существования искомых поднаборов точек утверждался более сильный факт о нечетности числа искомых поднаборов. Предложенное нами доказательство обобщения теорем~\ref{CGS} и~\ref{vKF} может быть обобщено на случай таких более общих формулировок.
\end{remark*}

Сформулируем теперь обобщение теорем~\ref{CGS} и~\ref{vKF}.

Будем говорить, что точки $A_1,A_2,\ldots,A_m\in\mathbb{R}^{n}$, $m>n$, находятся {\it в общем положении}, если никакие $n+1$ из них не лежат в одной гиперплоскости (то есть в аффинном подпространстве размерности $n-1$).

Пусть $\Delta$ и $\Delta'$ --- два $n$-мерных симплекса в $\mathbb{R}^{2n-1}$, вершины которых находятся в общем положении. Тогда эти симплексы {\it зацеплены}, если граница симплекса $\Delta$ пересекает внутренность симплекса $\Delta'$ в единственной точке.

\begin{theorem}
Пусть даны $n+3$ точки общего положения в $\mathbb{R}^{n}$. Тогда, если $n$ нечетно, то существуют два непересекающихся поднабора по $(n+3)/2$ точек, внутренности выпуклых оболочек которых пересекаются. Если же $n$ нечетно, то существуют два непересекающихся поднабора по $(n+2)/2$ точек, выпуклые оболочки которых зацеплены.
\label{main}
\end{theorem}

\section{Доказательство теоремы~\ref{main}}
\label{proof_main}

В доказательстве теоремы мы будем пользоваться следующим эквивалентным определением общности положения точек.

Точки $A_1,A_2,\ldots,A_m\in\mathbb{R}^{n}$, $m>n$, находятся {\it в общем положении}, если для любых $B_0,B_1,\ldots,B_{n}$ из них векторы $B_1-B_0,B_2-B_0,\ldots,B_{n}-B_{0}$ линейно независимы.

\subsection{Начало доказательства}
\label{proof_begin}
Обозначим точки из условия теоремы через $A_1,A_2,\ldots,A_{n+3}$. Рассмотрим следующую систему линейных уравнений относительно вещественных переменных $x_1,x_2,\ldots,x_{n+3}$:
\begin{equation}
 \begin{cases}
   x_{1}A_{1}+x_{2}A_{2}+\ldots+x_{n+3}A_{n+3}=0,
   \\
   x_{1}+x_{2}+\ldots+x_{n+3}=0.
 \end{cases}
\label{system}
\end{equation}

Поскольку в этой системе $n+1$ скалярных уравнений и $n+3$ неизвестных, то размерность пространства ее решений не меньше двух. А в силу общности положения точек $A_1,A_2,\ldots,A_{n+3}$ любые $n+1$ столбцов матрицы данной системы линейно независимы. Значит, размерность пространства решений в точности равна двум. Таким образом, данная система уравнений задает двумерную плоскость в $(n+3)$-мерном пространстве параметров. Обозначим эту плоскость через $\gamma$. Покольку любые $n+1$ столбцов матрицы системы~\eqref{system} линейно независимы, то при добавлении к системе~\eqref{system} любого из $n+3$ уравнений вида
$$
  x_{i}=0,
$$
где $i\in\{1,2,\ldots,n+3\}$, размерность пространства решений уменьшится на 1. Иначе говоря, пространство решений дополненной системы --- прямая в плоскости $\gamma$, проходящая через точку $(0,0,\ldots,0)$. Обозначим эту прямую через $\ell_i$. В силу линейной независимости любых $n+1$ столбцов матрицы системы~\eqref{system} прямые $\ell_1,\ell_2,\ldots,\ell_{n+3}$ попарно различны.

Рассмотрим в плоскости $\gamma$ произвольную точку $x=(x_1,x_2,\ldots,x_{n+3})$, не лежащую ни на одной из прямых $\ell_1,\ell_2,\ldots,\ell_{n+3}$. Обозначим через $\omega$ дугу окружности в плоскости $\gamma$ с центром в точке 0, концами которой являются точки $x$ и $-x=(-x_1,-x_2,\ldots,-x_{n+3})$. Для любого $i=1,2,\ldots,n+3$ прямая $\ell_i$ разделяет плоскость $\gamma$ на две полуплоскости, в одной из которых $i$-ые координаты точек положительны, а в другой --- отрицательны. Следовательно, при движении по дуге $\omega$ от точки $x$ к точке $-x$ при переходе через любую из прямых $\ell_1,\ell_2,\ldots,\ell_{n+3}$ число положительных координат точки изменяется ровно на 1.

\begin{lemma}
\label{main_lemma}
Пусть ненулевая точка $\hat{x}=(\hat{x}_1,\hat{x}_2,\ldots,\hat{x}_{n+3})$ является решением системы~\eqref{system}. Тогда выпуклые оболочки множеств $\{A_{i},\: \hat{x}_{i}>0\}$ и $\{A_{i},\: \hat{x}_{i}<0\}$ пересекаются по внутренней для обеих выпуклых оболочек точке.
\end{lemma}

\begin{proof}[Доказательство леммы.]
Обозначим через $I^{+}$ множество таких индексов $i\in\{1,\dots,{n+3}\}$, 
для которых $\hat{x}_{i}>0$, через $I^{-}$ --- множество индексов $i\in\{1,2,\ldots,n+3\}$, для которых $\hat{x}_{i}<0$. В силу того, что $\hat{x}$ является решением системы~\eqref{system}, выполнены равенства
\begin{equation*}
 \begin{cases}
   \sum\limits_{i\in I^{+}} \hat{x}_{i}A_{i} = \sum\limits_{i\in I^{-}} (-\hat{x}_{i}A_{i}),
   \\
   \sum\limits_{i\in I^{+}} \hat{x}_{i} = \sum\limits_{i\in I^{-}} (-\hat{x}_{i})\eqdef S.
 \end{cases}
\end{equation*}

Следовательно,
\begin{equation*}
 \sum\limits_{i\in I^{+}}\frac{\hat{x}_{i}}{S} A_{i} =
 \sum\limits_{i\in I^{-}}\frac{-\hat{x}_{i}}{S} A_{i}.
\end{equation*}

Это означает, что выпуклые оболочки наборов точек $\{A_{i},\: i\in I^{+}\}$ и $\{A_{i},\: i\in I^{-}\}$ пересекаются. Поскольку все числа $x_{i}$, $i\in I^{+}\cup I^{-}$, ненулевые, то точка пересечения является внутренней для обеих выпуклых оболочек.
\end{proof}

\subsection{Случай четного $n$}
\label{proof_even}
Поскольку у одной из точек $x$ и $-x$ положительных координат не меньше, чем $(n+2)/2$, а у другой --- не больше, чем $(n+2)/2$, то существует такое $i\in\{1,2,\ldots,n+3\}$, что у точки пересечения дуги $\omega$ с прямой $\ell_i$ будет ровно по $(n+2)/2$ положительных и отрицательных координат. Обозначим эту точку пересечения за $\hat{x}$. Тогда по лемме~\ref{main_lemma} выпуклые оболочки наборов $\{A_{i},\: \hat{x}_{i}>0\}$ и $\{A_{i},\: \hat{x}_{i}<0\}$, являющиеся $(n+2)/2$-мерными симплесками, пересекаются по внутренней точке.

\subsection{Случай нечетного $n$}
\label{proof_odd}
Очевидно, что два $(n+1)/2$-мерных симплекса $\Delta$ и $\Delta'$ в $\mathbb{R}^{n}$, вершины которых находятся в общем положении, зацеплены тогда и только тогда, когда они пересекаются по отрезку, причем один из концов этого отрезка является внутренней точкой для симплекса $\Delta$ и граничной точкой для симплекса $\Delta'$, а другой конец, наоборот, --- внутренней точкой для симплекса $\Delta'$ и граничной точкой для симплекса $\Delta$.

Поскольку у одной из точек $x$ и $-x$ положительных координат не меньше, чем $(n+3)/2$, а у другой --- не больше, чем $(n+3)/2$, то существуют такие две точки пересечения дуги $\omega$ с прямыми $\ell_1,\ell_2,\ldots,\ell_{n+3}$, что между ними на дуге $\omega$ нет других точек пересечения с прямыми $\ell_i$, и, кроме того, у одной из этих точек $(n+3)/2$ положительных и $(n+1)/2$ отрицательных координат, а у другой $(n+3)/2$ отрицательных и $(n+1)/2$ положительных координат. Обозначим эти точки соответственно через $x^{1}$ и $x^{2}$. Обозначим через $\Delta_{j}^{+}$, $j\in\{1,2\}$, выпуклую оболочку множества таких точек $A_i$, что $x_i^{j}$ положительно, а через $\Delta_{j}^{-}$, $j\in\{1,2\}$, --- выпуклую оболочку множества таких точек $A_i$, что $x_i^{j}$ отрицательно. В силу определения точек $x^{1}$ и $x^{2}$ выпуклые оболочки $\Delta_{1}^{+}$ и $\Delta_{2}^{-}$ являются $(n+1)/2$-мерными симплексами, причем $\Delta_{1}^{-}$ является $(n-1)/2$-мерной гранью симплекса $\Delta_{2}^{-}$, а $\Delta_{2}^{+}$ является $(n-1)/2$-мерной гранью симплекса $\Delta_{1}^{+}$. По лемме~\ref{main_lemma} существует внутренняя точка симплекса $\Delta_{1}^{+}$, принадлежащая симплексу $\Delta_{1}^{-}$, и существует внутренняя точка симплекса $\Delta_{2}^{-}$, принадлежащая симплексу $\Delta_{2}^{+}$. Что и означает зацепленность симплексов $\Delta_{1}^{+}$ и $\Delta_{2}^{-}$.

\end{document}